\font\steptwo=cmb10 scaled\magstep2
\font\stepthree=cmb10 scaled\magstep4
\magnification=\magstep1
\settabs 18 \columns
 
\hsize=16truecm
\baselineskip=17 pt

\def\s{\sigma}

\def\b{\bigskip}
\def\bb{\bigskip\bigskip}

\def\no{\noindent}

\def\ce{\centerline}
\def\ve{\vfill\eject}

  \font\got=eufm8 scaled\magstep1
 \def\g{{\got g}}


\def\Cit{\hbox{\it l\hskip -5.5pt C\/}}

\def\Crm{\hskip0.5mm \hbox{\rm l\hskip -5.5pt C\/}}

{\ce {\stepthree   QUASI HOPF DEFORMATIONS}}
\b {\ce{\stepthree OF QUANTUM GROUPS}}
 
\bb {\ce {C. Fr\o nsdal}}
 \b {\ce {Physics Department, University of California, Los Angeles CA 90095-1547, USA}}
\bb\bb\bb

\ce{{\it ABSTRACT}}

The search for elliptic quantum groups leads to a modified quantum Yang-Baxter relation and to a special class of
quasi-triangular quasi Hopf algebras. This paper calculates deformations of standard quantum groups (with or without
spectral parameter) in the category of quasi-Hopf algebras. An earlier investigation of the   deformations of
quantum groups, in the category of Hopf algebras, showed that quantum groups are generically rigid: Hopf algebra
deformations exist only under some restrictions on the parameters. In particular,  affine Kac-Moody algebras  are
more rigid than their loop algebra quotients;   and only the latter (in the case of sl(n)) can be deformed to
elliptic Hopf algebras. The generalization to quasi-Hopf deformations lifts this restriction. The full elliptic
quantum groups (with central extension) associated with $sl(n)$ are thus quasi-Hopf algebras. The universal
R-matrices satisfy a modified Yang-Baxter relation and are calculated more or less explicitly.  The modified
classical Yang-Baxter relation  is obtained, and the elliptic solutions are worked out explicitly.

The same method is used to construct the Universal R-matrices associated with Felder's  quantization of the
Knizhnik-Zamolodchikov-Bernard equation, to throw some light on the quasi Hopf structure of conformal field theory
and (perhaps)    the Calogero-Moser models.

\def\Rt{\tilde R}

\ve

  {\steptwo 1. Introduction.}
\b The pioneering work of J.R. Baxter [Ba] on the 8-vertex model drew attention to the fundamental importance of the
elliptic R-matrix associated with $sl(2)$. The 6-vertex model is a special limiting case, in which the coefficients
of the R-matrix are trigonometric functions of the spectral parameter. The generalization to $sl(n)$ was found by
Belavin [Be]. 

The work of  Drin'feld [D1] and many others showed that the trigonometric R-matrix has a beautiful interpretation in
terms of   quantized affine Kac-Moody algebras, its role being to turn these algebras into Hopf algebras (quantum
groups).

For the physical applications, the distinction between the full (quantized) Kac-Moody algebra and the underlying
loop algebra is fundamental. The former is a central extension of the latter and this central extension is
essential; it lies behind the concept of ``level" of physical representations.

In the case of the elliptic R-matrix our understanding of it in terms of quantized Kac-Moody algebras
 is, so far, limited to the loop algebra quotient; that is, to level zero.  The $RLL$ formalism defines a dual
algebra of physical observables by
$$ R_{12}L_1L_2 = L_1L_2R_{12},\eqno(1.1)
$$ but this method is applicable only at level zero, since the R-matrix is known only in finite dimensional
representations that annihilate the center, and since the universal R-matrix is known only for the quotient algebra.
Attempts have been made [Fo] to define a full elliptic quantum group with the help of  a modified relation, of the
type
$$ R_{12}L_1L_2 = L_1L_2R_{12}^*.\eqno(1.2)
$$ This produces an algebra but,  as far as is known, not a co-product, which hints at a connection to modified
Yang-Baxter relations and quasi-Hopf algebras that have appeared in another context [BBB]. In fact, such a
connection is established in this paper.

Recently, in a study  of deformations of ``standard" quantum groups [Fr], it was found, very unexpectedly, that the
universal elliptic R-matrix for $sl(n)$ is a deformation (by twisting) of the trigonometric one; however, this is
true only at  the level of loop algebras, without the central extension. In fact, the standard quantum groups are
generically rigid; deformations exist only for special values of the parameters, and  the elliptic deformation only
if the generator of the central extension is replaced by zero. 

In this paper it is reported that the generic rigidity of standard quantum groups is very much relaxed when
deformations into a wider category of quasi Hopf algebras is considered. The deformed R-matrix
$\tilde R$ satisfies a modified Yang-Baxter relation
$$
\Rt_{12}\Rt_{13,2}\Rt_{23} = \Rt_{23,1}\Rt_{13}\Rt_{12,3},\eqno(1.3)
$$ where $\Rt_{12} \in {\cal A}' \otimes {\cal A}'$ and $\Rt_{12,3} \in {\cal A}' \otimes {\cal A}' 
  \otimes {\cal A}_0'$. (${\cal A}'$ is the quantized Kac-Moody algebra and $ {\cal A}_0'$ is the Cartan
subalgebra.) This type of modified Yang-Baxter relation first \break appeared in the works of Gervais and Neveu [GN]
and Felder [Fe]; a quasi Hopf interpretation was given, in a   special case, by Babelon, Bernard and Billey [BBB].

This deformed R-matrix is related to the standard R-matrix by a twist, and the twisted coproduct gives  to the
algebra ${\cal A}'$ the structure of a bialgebra (non-coassociative in general), more precisely a quasi-triangular,
quasi-Hopf algebra. A possibility for a modified $RLL$ formalism is to postulate
$$
\Rt_{12}L_{13,2}L_{23} = L_{23,1}L_{13}\Rt_{12,3},\eqno(1.4)
$$  where now the third space is identified with the quantum space. In the most interesting elliptic case the
extension of $\Rt_{12}$ to the third space is supported on the center.  The algebra defined by (1.4) is associative;
therefore, it is not the quasi Hopf dual of ${\cal A}'$. The discovery that the  algebra of physical observables, in
the quasi Hopf case, is unrelated to the dual of the quantized Kac-Moody  algebra is an important turning point in
the development of the theory. 

This paper calculates the universal R-matrices for quasi-Hopf deformations of standard quantum groups, more or less
explicitly. The classical limit is taken, and the expressions obtained for the classical  r-matrices are quite
explicit.  Examples, besides the elliptic quantum groups (in the sense of the opening of this introduction), include
the ``elliptic quantum groups" of Felder [Fe], with  application to Calogero-Moser models and to the
Knizhnik-Zamolodchikov-Bernard equation [ABB][FR][Fe]. Finally, it may be  noted that the integrable models that
have a spectral parameter on a curve of higher genus are also governed by Eq.(1.4) [Au]. Perhaps it will turn out
that these models too have a quasi Hopf interpretation.
\ve

\no {\steptwo 2. Standard Quantum Groups and   Standard \break R-matrix.} (Review.)
\b This section reviews, as brifly as possible, the first part of [Fr], to fix some notations and explain the point
of view that we continue to adopt here. First of all, the algebras.
\b {\bf Definition 2.1.}   Let $M,N$ be two countable sets,
$\varphi,\psi$  two maps,
$$
\eqalign{& \varphi :~~M\times M \rightarrow \Crm~, \cr & \psi :~~M\times N \rightarrow \Crm~, \cr} \quad
\eqalign{a,b &\mapsto \varphi^{ab}~, \cr a,\beta & \mapsto H_a(\beta)~. \cr} \eqno(2.1)
$$
\no Let ${\cal{A}}$ or ${\cal{A}}(\varphi,\psi)$ be the universal, associative, unital  algebra over \Crm \enskip
with generators
$\{H_a\}\, a\in M,~ \{e_{\pm \alpha}\}\,\alpha \in N$, and relations
$$
\eqalignno{&[H_a,H_b]=0~, \quad  [H_a,e_{\pm\beta}] = \pm H_a(\beta)e_{\pm\beta}~, & (2.2) \cr
&[e_\alpha,e_{-\beta}]=\delta^\beta_\alpha
\bigl(e^{\varphi(\alpha,\cdot)}-e^{-\varphi(\cdot,\alpha)}\bigr)~, &  (2.3) \cr}
$$
\no with $\varphi(\alpha,\cdot)=\varphi^{ab}H_a(\alpha)H_b,~
\varphi(\cdot,\alpha)=\varphi^{ab}H_aH_b(\alpha)$ and $ e^{\varphi(\alpha, \cdot) + \varphi(\cdot,\alpha)} \neq 1,
~~ \alpha \in N
$.
\bb  Our actual concern in this paper is with  quotient algebras ${\cal A}'$ that will be introduced presently. But
first, the standard R-matrix.
\b

\no {\bf Definition 2.2.}  A standard R-matrix is a formal series of the form
$$
\eqalign{R &= \exp\bigl(\varphi^{ab}H_a\otimes H_b\bigr)
\bigl(1+e_{-\alpha}\otimes e_\alpha  + \sum^\infty_{k=2} t^{\alpha^\prime_1\ldots a^\prime_k}_ {\alpha_1\ldots
\alpha_k} e_{-\alpha_1}\ldots e_{-\alpha_k}\otimes e_{\alpha^\prime_1}\ldots e_{\alpha^\prime_k}\bigr)~. \cr}
\eqno(2.4)
$$
\b
\no In this formula, and in others to follow, summation over repeated indices is implied. 
  For fixed
$(\alpha)=\alpha_1,\ldots,\alpha_k$ the sum over $(\alpha^\prime)$ runs over the permutations of $(\alpha)$. The
coefficients
$t^{(\alpha^\prime)}_{(\alpha)}$ are in \Crm. 
\b When the parameters of the algebra ${\cal A}$ are in general position, then there is a unique, standard R-matrix
that satisfies the Yang-Baxter relation. That is, this relation fixes all the  coefficients in $R$, and uniquely.
This was proved [Fr] by showing that the Yang-Baxter relation for $R$ is  equivalent to a recursion relation for the
coefficients, and that the said recursion relation is generically uniquely solvable. For special values of the
parameters this is not the case, but then the obstruction is avoided by passing to the quotient ${\cal A}'$ by a
certain ideal   ${\cal I}\in {\cal A}$. In the case of quantized Kac-Moody algebras the ideal is generated by the
Serre-Drin'feld relations. The details need not be reviewed here; we are mainly interested in Kac-Moody algebras and
pass directly to the standard R-matrix on ${\cal A}'$, on which (2.4) remains valid.
\bb
\bb

\no {\steptwo 3. First Order Deformations to Quasi Hopf.}
\b A Hopf deformation of the standard R-matrix is a formal series
$$ 
\Rt = R_\epsilon=R+\epsilon R_1+\epsilon^2 R_2 + \ldots~. \eqno(3.1)
$$
\no Here $R$ is a standard R-matrix on ${\cal A}' = {\cal A}/{\cal I}$ with any choice of  parameters and the ideal
${\cal I}$ determined by them. The coefficients  $t_{(\alpha)}^{(\alpha')}$ of $R$ are determined by the Yang-Baxter
relation, and we attempt to find
$R_1,R_2,\ldots$ so that $R_\epsilon$ will satisfy the same relation to each order in
$\epsilon$. 

To describe the deformations calculated in [Fr] we need to introduce a grading in ${\cal A}' \otimes {\cal A}'
\otimes{\cal A}'$. First, ${\cal A}'$ admits a grading in which the generators $e_\alpha$ have grade $1$, the
$e_{-\alpha}$ have grade $-1$ and the $H_a$ have grade zero. This induces a natural grading in ${\cal A}'\otimes
{\cal A}'$.    Each term in $R$ has a grade $(a,b)$ defined as grade$(a,b) = (-$grade$(a), ~$grade$(b))$, with $a,b
\geq 0$, the lowest grade is
$(0,0)$. The lowest grade in 
$R_\epsilon$ and in $R_1$ is $(-1,-1)$. In fact, the deformation is completely fixed by the term of lowest grade in
$R_1$, we call it the driving term, and the Yang-Baxter relation. To grade the elements $(a,b,c) \in {\cal A}'
\otimes {\cal A}'
\otimes{\cal A}'$ we use  
$$ {\rm grade} (a,b,c) := (-{\rm grade}(a),~ {\rm grade}(c)).
$$  

The deformations    are elementary or compound. An elementary deformation has a driving term of the type
$$
 S~e_\s \otimes e_{-\rho}, \quad S \in {\cal A}'_0 \otimes {\cal A}'_0. \eqno(3.2)
$$
\b

\no {\bf Theorem 3.1.} [Fr]  Let $R$ be the standard R-matrix for ${\cal A}'$.
  Suppose that $R+\epsilon R_1$ is a first order deformation, satisfying the Yang-Baxter relation to first order in
$\epsilon$.  Suppose also that the term of lowest grade in $R_1$ has the form (3.2); then  the parameters satisfy
$$  e^{\varphi(\cdot,\rho)+\varphi(\sigma,\cdot)}=1~. \eqno(3.3)
$$  Conversely, when the parameters satisfy (3.3) for some pair $(\s,\rho )$, then there exists a unique first 
order deformation such that the term of lowest grade has the form (3.2), namely
$$  R_1=  R(Ke_\sigma\otimes K'e_{-\rho}) - (K'e_{-\rho}\otimes Ke_\sigma)R ~, \eqno(3.4)
$$
\no with $K = e^{-\varphi(\s,.)},~K'~:=~e^{\varphi(\cdot,\rho)}$.
\b  The first order deformations form a linear space. The elementary deformations form a basis. In the case of
compound deformations there are obstructions in higher orders of the deformation parameter $\epsilon$ to which we
shall turn in Section 4. 

The necessity of the condition (3.3) is established easily, by an examination of the terms of  lowest grade in the
Yang-Baxter relation. One encounters the obstruction
$$ 1 - Q(\s,\rho),  \quad Q(\s,\rho) := e^{-\varphi(\s,.)-\varphi(.,\rho)}\eqno(3.5) 
$$ in the second space. That is; standard quantum groups are rigid (as quantum groups) to deformations, unless the
condition (3.3) is satisfied for one or more pairs $(\s,\rho)$.  

The window defined by the condition (3.3) is wide enough to admit all known simple quantum groups. There is,
nevertheless, a strong  motivation for relaxing this condition. Consider the case of a quantized, affine Kac-Moody
algebra with   generators $c,d$ of central extension and scale, and imaginary root generators
$e_0,e_{-0}$. Then  
$$
\varphi(.,.) = \sum \varphi^{ab}H_a \otimes H_b + u\, c\otimes d + (1-u)\, d \otimes c,\eqno (3.6)
$$ with the parameter $u\in \Cit$.  Choose a pair 
$(e_\s,e_{-\rho})$ of the type
$$ e_\s = e_0, \quad e_{-\rho} \neq e_{-0},\eqno(3.7)
$$  then, since $d(\s) = \delta_\s^0$,
$$
\varphi(\s,.) = \sum \varphi^{ab}H_a(\s) \otimes H_b + (1-u)\, c,
\quad
\varphi(.,\rho) = \sum \varphi^{ab}H_a \otimes H_b(\rho).
$$  Evidently, (3.3) leads to $u=1$, which, by the way, demonstrates the non-trivial nature of this parameter, first
included in [Fr]. But the deformations that lead to elliptic R-matrices are compound deformations that include pairs
of the type (3.7) as well as the opposite type,
$ e_\s \neq e_0, ~ e_{-\rho} = e_{-0},
$ and a similar calculation now leads to $u = 0$. The options are (1) pass to the quotient defined by setting $c=0$
or (2) relax the condition (3.3). Here, of course, we choose  to relax the condition on the parameters.

As we saw, this implies relaxing the Yang-Baxter relation. In fact, it is not difficult to see that the 
modification that is required is of the type
$$
\Rt_{12}\Rt_{13,2}\Rt_{23} = \Rt_{23,1}\Rt_{13}\Rt_{12,3},\eqno(3.8)
$$ with $R_1$ of the type (3.4) and
$$ R_{12,3} = R_{12}\otimes 1, \quad (R_1)_{12,3} = (R_1)_{12} \otimes Q(\s,\rho).\eqno(3.9)
$$
\b {\bf Theorem 3.2.} (a) Let $R$ be the standard R-matrix,  then the deformed R-matrix \break $R + \epsilon R_1$,
with $R_1$ as in (3.4) and $R_{12,3}$   as in (3.9), satisfies the modified Yang-Baxter relation (3.8), to first
order in
$\epsilon$. (b) If the term of lowest grade in $R_1$ has the form (3.2),   then both
$R_1$ and
$(R_1)_{12,3}$ are determined uniquely by (3.8).
\b A direct proof follows the plan of the proof of Theorem 3.1, in [Fr].   A simpler proof of (a) is obtained in 
Section 4; therefore the important fact to notice here is part (b), since it shows that all these first order
deformations can be interpreted as twists in the sense of [D3]. This justifies the term  quasi Hopf in the heading
of this section.
 
\ve

\no {\steptwo 4. Formal Deformations.}
\b Here we study deformations by formal power series, to all orders in the deformation parameter
$\epsilon$, sometimes referred to as exact deformations since a first order deformation may be thought of as an
approximation. However, not all first order deformations can be continued to higher powers of
$\epsilon$.  We have seen that first order deformations are twists, and we shall now construct formal deformations
of the same kind.

 For the following result ${\cal A}'$ is any coboundary Hopf algebra.
\b
\no {\bf Theorem 4.1.} Let $R$ be the R-matrix, and $ \Delta$  the coproduct, of a coboundary Hopf algebra ${\cal
A}',$ and $F \in {\cal A}'\otimes {\cal A}'$, invertible, such that
$$
\bigl((1 \otimes \Delta_{21}) F\bigr)F_{12} = \bigl( (\Delta_{13} \otimes 1) F\bigr) F_{31};
\eqno(4.1)
$$
\vskip-.2cm
\no then
\vskip-.3cm
$$
\tilde R := (F^t)^{-1} R F\eqno(4.2)
$$ (a) satisfies the Yang-Baxter relation and (b) defines a coboundary Hopf algebra $\tilde {\cal A}$ with the same
product and with co-product
$$
\tilde \Delta = (F^t)^{-1} \Delta F^t.\eqno(4.3)
$$
\b
 This result is due to Drin'feld; a simple proof was given in [Fr].
 The cocycle condition (4.1) first appears in [G].
 
Without  (4.1), one can define an invertible element $\Phi \in {\cal A}'^{\otimes 3}$ by the formula
$$
\bigl((1 \otimes \Delta_{21}) F\bigr)F_{12} = \bigl( (\Delta_{13} \otimes 1) F\bigr) F_{31}\Phi.
\eqno(4.4)
$$  The coproduct (4.3) is then not in general coassociative, instead we have
$$
  (1 \otimes \tilde \Delta )   \tilde \Delta (a) \Phi = 
\Phi (\tilde \Delta  \otimes 1\bigr) \tilde\Delta (a).\eqno(4.5)
$$ If we define
$$ F_{12,3} := F_{12}\Phi, \quad \tilde R := (F^t)_{-1}RF,\quad \tilde R_{12,3} := 
F_{21,3}^{-1}R_{12}F_{12,3}\eqno(4.6)
$$ then the modified Yang-Baxter relation (3.8) is satisfied by $\tilde R$.

If one does not impose a condition such as     (4.1), then one runs the risk of triviality.  Eq.s (4.2) and (4.3)
define a coboundary quasi-Hopf algebra. If $F$ is a formal power series with constant term 1, then $\Rt$ is a formal
power series with constant term $R$, and thus automatically a formal deformation of $R$. There is thus no
problematics.

The problem becomes much less trivial if one puts   conditions on the intertwiner.   The loss of associativity is
not welcome and it is natural to try to contain the damage.
\b
\no {\bf Definition 4.} We shall say that the twist $F$ is benevolent, resp. central, if $F$ and
$F_{12,3}$ have a  representation
$$ F = \sum_i F^i \otimes F_i, \quad F_{12,3} = \sum_i F^i \otimes F_i \otimes Q(i),\eqno(4.7)
$$ with $Q(i) \in {\cal A}'$, resp.  $Q(i)\in $ the center of ${\cal A}'$.
\b
 The twists that permit to define elliptic quantum groups (Section 5) are central.  The first order deformations of
Section 3 are all obtained by twists that are benevolent to first order, with the additional simplification that
$Q(i) \in {\cal A}'_0$. We shall construct formal deformations with the same property. 

It   is convenient to replace the generators $e_{\pm \alpha} $ by
$$ f_\s := e^{-\varphi(\s,.)}e_\s, \quad f_\rho := e_{-\rho}\,e^{\varphi(.,\rho)},\eqno(4.8)
$$ satisfying
$$ [f_\s, f_{-\rho}] = \delta_\s^\rho (e^{\varphi(.,\s)} - e^{-\varphi (\s,.)}) \eqno(4.9)  
$$ and
$$
\Delta f_\s = e^{-\varphi(\s,.)} \otimes f_\s + f_\s \otimes 1, \quad
\Delta f_{-\rho} = 1 \otimes f_{-\rho} + f_{-\rho} \otimes e^{\varphi(.,\rho)}.
$$ In the case of an elementary first order deformation (with $\s \neq \rho$),  it is easy to construct the
corresponding formal deformation. In this case $F$ and $\Phi$  take the form
$$ F = 1 - \epsilon ~f_\s \otimes f_{-\rho}\,Q + o(\epsilon^2), \quad Q = e^{-\varphi(\s,.)
-\varphi(.,\rho)},\eqno(4.10)
$$
$$ F\Phi = 1 - \epsilon~f_\s \otimes f_{-\rho}Q \otimes Q + o(\epsilon^2).\eqno(4.11)
$$
\b
\no {\bf Proposition 4.1.}  If $F$ is a formal, benevolent twist that agrees with (4.10) to first order in
$\epsilon$, and if $\s \neq \rho$, then (a)  
$$ Qf_\s f_{-\rho} = f_\s f_{-\rho}Q,\eqno(4.12)
$$ (b)
\vskip-10mm
$$ F = e_q^{-\epsilon f_\s \otimes \,f_{-\rho}Q} = \sum {1 \over [n!]_q}(-\epsilon f_\s \otimes
f_{-\rho}Q)^n,\eqno(4.13)
$$  and
$$ F\Phi =  e_q^{-\epsilon f_\s \otimes  \,f_{-\rho}Q \,\otimes \,Q},\quad q = e^{-\varphi(\s,\s)}.\eqno(4.14)
$$
\b
\no {\bf Proof.} (a) That the condition (4.12) is necessary will be shown in greater generality below. (b) Using
(4.13) we have
$$
\bigl((1 \otimes \Delta_{21})F\bigr)F_{12} = e_q^{(A+B)}e_q^C,
\quad \bigl((\Delta_{13} \otimes 1)F\bigr)F_{31,2} = e_q^{B+C}e_q^A,
$$ with $A = -\epsilon f_{-\rho}Q \otimes Q \otimes f_\s,~ B = -\epsilon e^{-\varphi(\s,.)} \otimes f_{-\rho}Q
\otimes f_\s, ~ C = -\epsilon f_\s  \otimes f_{-\rho}Q \otimes 1$.  Here  $BA = qAB,  ~CB = qBC$ and  $AC = CA$ 
because $\s \neq \rho$; both expressions reduce to $e_q^Be_q^Ae_q^C$. 
\b The need to exclude the possibility of $\s = \rho$ arises in the quasi-Hopf case; in the case of a Hopf
deformation, when $Q = 1$, this possibility is excluded {\it a priori} by the last condition in Definition (2.1).
The   case when $\sigma = \rho$ will be solved below, Eq.(4.30).

The general case of compound deformations is much more difficult. First of all, we encounter  additional
obstructions in higher orders of $\epsilon$.
\b
\no {\bf Proposition 4.2}. (Initial Conditions.) Let $[\tau]$ be a set of pairs $(\s,\rho) \in N \times N$, and $F$ a
formal, benevolent twist of the type, \footnote *{Here and below, the parameters $\epsilon_\s$ should be understood
to stand for $k_\s\epsilon$, with fixed parameters $k_\s$ and a single deformation parameter $\epsilon$.}
$$ F = 1 -\sum_{(\s,\rho) \in [\tau]}\hskip-1mm\epsilon_\s ~f_\s \otimes f_{-\rho}Q(\s,\rho) + o(\epsilon^2),\quad
Q(\s,\rho) = e^{-\varphi(\s,.) -\varphi(,.\rho)}.\eqno(4.15) 
$$ Then (a) the set $[\tau]$ is the restriction to simple roots of the graph of an isomorphism 
$\tau: \Gamma_1 \mapsto \Gamma_2$, where $\Gamma_1, \Gamma_2$ are subalgebras of ${\cal A}'$ generated by simple
roots, and (b)
$$ Q(\s',\rho')f_\s f_{-\rho} = f_\s f_{-\rho}Q(\s',\rho'),\quad (\s,\rho), (\s',\rho') \in [\tau].\eqno(4.16)
$$
\b The proof of (a) is an easy adaptation of Proposition 15.1 of [Fr].  It consists of two steps. First, a recursion
relation for $F$ is derived from the modified Yang-Baxter relation. Then the  integrability of this recursion
relation is shown to require the property of the set $[\tau]$ stated in the theorem. In the
 most interesting cases (elliptic quantum groups),  $\Gamma_1 = \Gamma_2$ is generated by all the simple roots. (b)
The necessity of this condition will be proved in the course of the proof of the  next theorem. (The statement is
included here so that all conditions required by the term of first order in $\epsilon$ are collected in one place.)

We now turn our attention to the explicit determination of $F$,  solving the cocycle condition (4.4) in terms of
formal series with the initial condition (4.15).

From now on we shall always restrict the parameters so that, if $(\s,\rho) \in [\tau]$, and $\rho \in  
\Gamma_1$, then $\epsilon_\s =
\epsilon_\rho$. This can be arranged by a renormalization of the generators.
 \b
\no {\bf Theorem 4.2.} (a) If the formal series $F\in {\cal A}'
\otimes {\cal A}'$ and  $F\Phi\in {\cal A}'
\otimes {\cal A}' \otimes {\cal A}'$ satisfiy the cocycle condition (4.4), as well as the initial conditions (a) and
(b) of Proposition 4.2, then they have the form
$$ F = \sum_{\s,\rho}[\s]\otimes [-\rho]T(\s,\rho),\quad F\Phi = \sum_{\s,\rho} [\s] \otimes [-\rho] T(\s,\rho)
\otimes Q(\s,\rho).\eqno(4.17) 
$$ Here and below  we use a multi-index notation; $\s$~($\rho$) is an abbreviation for a  set of elements taken from
$\Gamma_1$ ($\Gamma_2$), $(\s,\rho)$ is a set of the form $(\{\s_1,\s_2,...\}, 
\{\tau^{m_1}\s_1,\tau^{m_2}\s_2,...\})$ for some positive integers numbers $m_i$ such that $\tau^{m_i} \s_i \in
\Gamma_2$. In $[\s] \otimes [-\rho]$,   $[\s]$ is a product of  $f_\s$'s and $ [-\rho]$ is a product of    
$f_{-\rho}$'s and a numerical factor. The numerical factors are uniquely determined by the following  recursion
relation,
$$
\bigl(1 \otimes e^{\varphi(,.\rho)}\partial_\rho\bigr)F + \sum_{\tau^m\s = \rho} \epsilon_\s^m\,  [ 1 \otimes f_\s,
F]_q + \sum_{\tau^m\s = \rho} \epsilon_\s^m\,\bigl( f_\s \otimes e^{-\varphi(\s,.)}\bigr)F = 0,
\eqno(4.18)
$$
 Here $ [A,B]_q = AB - q BA$ and $q \in \Cit$ is determined by
$Q f_\s = qf_\s Q $.
 (b) The factor $Q(\s,\rho)\in {\cal A }'_0$  is fixed by two properties: (1) when the  set $\s$ has only one
element, and $\rho = \tau \s$, then (3.5) applies, and (2) the groupoid property
$$ Q(\s,\rho\cup\s')Q( \s', \rho') = Q(\s,\rho\cup\rho').\eqno(4.19)
$$

\no (c) The factor $T(\s,\rho)\in {\cal A }'_0$ is
$$ T(\s,\rho) = -\prod_i\epsilon_{\s_i}^{m_i} \,Q(\s,\rho),\eqno(4.20)
$$ where the product runs over $(\s_i,m_i)$ in    $\rho =
\{\tau^{m_1}\s_1,\tau^{m_2}\s_2,...\}$.  
\b   
 \no {\bf Proof.} The cocycle condition reads
$$
\eqalign{\bigl(\sum\Delta'[-\rho]T({\s,\rho}) &\otimes [\s]\bigr)\bigl(\sum[\s] \otimes [-\rho]T(\s,\rho)
\otimes 1\bigr)\cr &= \bigl(\Delta_{13}[\s] \otimes [-\rho]T(\s,\rho)\bigr)
\bigl(\sum[-\rho]T(\s,\rho) \otimes Q(\s,\rho) \otimes [\s]\bigr).\cr}\eqno(4.21)
$$
   Balancing terms with no roots in the second space gives us 
$$
\Delta T = T \otimes Q.\eqno(4.22)
$$  Balancing terms with one root in the second space gives
$$
\eqalign{
\sum[-\rho]T(\s,\rho)f_{\s'} &\otimes Q(\s,\rho)f_{-\rho'}T(\s',\rho') \otimes [\s]\cr &+ \sum
e^{\varphi(\s',.)}\bigl(\partial_{\rho'}[-\rho]T(\s,\rho) \otimes f_{-\rho'}Q(\s,\rho) \otimes [\s]\cr 
   }
$$
\vskip-9mm
$$  \eqalign{ = \sum e^{-\varphi(\s',.)}[-\rho]T(\s,\rho)&\otimes f_{-\rho'}T(\s',\rho')Q(\s,\rho) \otimes
f_{\s'}[\s]\cr &+ \sum f_\s[-\rho]T(\s,\rho) \otimes f_{-\rho'}T(\s',\rho')Q(\s,\rho) \otimes [\s]. 
\cr}\eqno(4.23)
$$ The operator $e^{\varphi(\s',.)}\partial_{\rho'}$,  substitutes $e^{\varphi(\s',.)}$ for $f_{-\rho'}$ in
$[-\rho]$. In order that the first term on the left combine with the second term on the right we need
$$ Q(\s,\rho)f_\s f_{-\rho} =f_\s f_{-\rho}Q(\s,\rho).\eqno(4.24)
$$ This is a generalization of (4.12) and (4.16) and completes the proofs of Propositions 4.1 and 4.2.
 By inspection of the factors that appear in the second space
  we find that $Q$ and $T$ must have the groupoid property, which proves parts (b) and (c). Finally, (4.23) now
reduces to (4.18). 
 \b This recursion relation is (uniquely) solvable.  In fact, the solution is
$$ F = F^1F^2....,\eqno(4.25)
$$
\vskip-2mm 
\no with $F^m$ of the form
$$ F^m = \sum_{\sigma} [\s] \otimes [-\tau^m\s]T(\s,\tau^m\s) = 1 - \sum_\sigma \epsilon_\s^mf_\s \otimes
f_{-\tau^m\s} Q(\s,\tau^m\s) + o(\epsilon^{2m}),\eqno(4.26)
$$ determined by the recursion relation
$$ (1 \otimes e^{\varphi(.,\rho)}\partial_\rho) F^m + \epsilon_\s^m   \bigl(f_{\tau^{-m}\rho} \otimes
e^{\varphi(.,\rho)}Q(\tau^{-m}\rho,\rho)\bigr) F^m = 0.\eqno(4.27)
$$ The unique solution is the same as in the Hopf case, except for the $Q$-factors. 

A compact expression for $F\Phi$ is
$$ (F\Phi)_{123} = e^{\varphi_{32}-\varphi_{13}} F_{12}e^{\varphi_{13}-\varphi_{32}}.\eqno(4.28)
$$
\b
\ce {\it EXAMPLE}
 
Consider a deformation of $sl_q(2)$ (or any quantized Kac-Moody algebra), with
$$ F = 1 -\epsilon f_1 \otimes f_{-1}e^{-K} + o(\epsilon^2),\quad K =  \varphi(1,.) + \varphi(.,1).
\eqno(4.29)
$$ Then $F$ and $F^m$ have the form
$$ F = \sum_{i=0}^\infty f_1^i \otimes f_{-1}^i\Psi_i  = \prod_{m=1}^\infty  F^m, \quad F^m = \sum_{i=0}^\infty A_i^m
~f_1^i
\otimes f_{-1}^ie^{-imK},\eqno(4.30)
$$ with coefficients $\Psi_i \in {\cal A}'$ and $A_i^m \in \Cit$. The recursion   (4.27) gives, with  $q =
e^{\varphi(1,1)}$,
$$ A_i^m = {(-\epsilon^m)^i\over [i]_q} q^{i(i-1)m}, \quad F^m = e_q^{-\epsilon^mf_1 \otimes f_{-1}e^{-mK}}.
$$ We can also  compute $F$ directly using (4.18). A very short calculation leads to
$$
\Psi_i = {(-\epsilon)^i\over [i!]_q}{q^{i(i-1)} e^{-iK}\over \prod_{\alpha = 0}^{i-1}(1 - q^\alpha\epsilon e^{-K})}.
$$ This result was first reported  in [BBB]. Finally,
$$ F = e_q^{-\epsilon f_1 \otimes f_{-1}Q(1-\epsilon Q)^{-1}},\quad Q = e^{-K}.\eqno(4.31)
$$ 
\b The relation (4.27) is equivalent to the following:
$F^m$ is the unique solution (as a power series with the given boundary conditions (4.26)) of
$$
 [1 \otimes f_\rho, F^m]_q +   \epsilon^m \bigg\{ \bigl(f_{\tau^{-m}\rho} \otimes e^{\varphi(., \rho)}Q
 \bigr) F^m  -  qF^m\bigl( f_{\tau^{-m}\rho} \otimes e^{-\varphi( \rho,.)} Q \bigr)\bigg\} = 0. 
\eqno(4.32)
$$ This last equation is more convenient for calculating the classical limit. 
\ve
\no { \steptwo 5. Elliptic quasi Lie algebras.}
\b If we introduce the parameter $\hbar$ in the usual way, and expand,
$$ R = 1 + \hbar r + o(\hbar^2),  \eqno(5.1)
$$ then the Yang-Baxter relation for $R$ gives the classical Yang-Baxter relation for $r$. The leading order is
$\hbar^2$, but terms of that order in $R$ cancel out in the Yang-Baxter relation and therefore need not be
calculated. However, to get the same accuracy in the modified Yang Baxter relation we need to carry the expansion a
little further, as follows,  
$$
\tilde R = 1 + \hbar \tilde r + o(\hbar^2), \quad Q = 1 + \hbar Y + o(\hbar^2), 
$$ and schematically, with $\tilde r = r + \delta r$,
$$
\tilde R_{12.3} = 1 + \hbar\, r + \hbar \, \delta r_{12}(1 + \hbar Y_3) + o(\hbar^2).\eqno(5.2)
$$ The unwritten $\hbar^2$ terms cancel among themselves, and the resulting modified classical Yang-Baxter relation
is, schematically
$$ [\tilde r_{12},\tilde r_{13}]   + [\tilde r_{12},\tilde r_{23}]   + [\tilde r_{13},\tilde r_{23}] +
\delta r
\wedge Y
  = 0.\eqno(5.3)
$$ To explain the fourth term we need to make some calculations. 

First of all, to obtain $\tilde r$,  set  $F^m = 1 + \hbar X^m  + o(\hbar^2)$,  then from (4.32) we get the relation
$$ [1 \otimes f_\rho + \epsilon^m\,f_\s \otimes 1, X^m] = -\epsilon^m f_\s \otimes (\varphi +
\varphi^t)(\rho), \quad \tau^m\s = \rho,\eqno(5.4)
$$
 as in the Hopf case except that  the restriction on the parameters has been relaxed.  

Next, the formula (4.28) gives
$$ F_{12,3} = 1 + \hbar X_{12} + \hbar^2[\varphi_{32} - \varphi_{13}, X_{12}] + o(\hbar^2),
$$ and 
$$
\delta r \wedge Y = - \sum_\s \epsilon_\s {d\over d \epsilon_\s} \delta r \wedge Y_\s,\quad Y_\s :=\varphi(\s,.) +
\varphi(.,\tau \s).\eqno(5.6) 
$$ 
\ve
\ce {\it EXAMPLES}

(1) The simplest example is   an elementary deformation of $sl_q(3)$, with
$$ F = 1 - \epsilon\,f_1 \otimes f_{-2}\,e^{-\varphi(1,.) - \varphi(.,2)} + o(\epsilon^2).\eqno(5.7)
$$
$$ F_{12,3} = 1 -\epsilon \,f_1 \otimes f_{-2}\,e^{-\varphi(1,.) - \varphi(.,2)} \otimes e^{-\varphi(1,.) -
\varphi(.,2)} + o(\epsilon^2).\eqno(5.9)
 $$ Terms beyond the first order in $\epsilon$   are irrelevant as far as the classical limit is concerned. (This is
because $\tau^m  \Gamma_1\cap \Gamma_1$ is empty for $m \geq 1$.) In this case the classical limit is exactly
$$
\tilde r = r + \epsilon\,\delta r,\quad \delta r = e_{-2} \wedge e_1,\eqno(5.8)
$$ where $r$ is the standard r-matrix $$ r = \varphi + \sum_{i=1}^2 e_{-i} \otimes e_i.\eqno(5.9)
$$
 There is only one term in $\delta r \wedge Y$, with 
$Y = -\varphi(1,.) -
\varphi(.,2)$, and the validity of (5.3) is easily verified directly. The special case   $Y = 0$ is an esoteric Lie
bialgebra [CG][FG]. 
\b (2) A   more complicated example is the deformation (4.29)
$$ F = 1 -\epsilon f_1 \otimes f_{-1}e^{-K} + o(\epsilon^2),\quad K =  \varphi(1,.) + \varphi(.,1).
$$ In this case $\Gamma_1 = \Gamma_2$, $\tau^m f_1 \in \Gamma_1$ for all $m$ and the factorization of $F$ is an
infinite product. However,  only the first two terms in the expression (4.29) for $F^m$ contribute to the classical
limit; this is because $\Gamma_1$ is abelian. Consequently,
$$
\tilde r = r + \sum_m \epsilon^m ~e_{-1}\wedge e_1 \eqno(5.10) 
$$
\vskip-.5cm
\no and 
$$
\delta r \wedge Y = -\sum_m m\epsilon^m~e_{-1}\wedge e_1\wedge K,\quad K = \varphi(1,.) +
\varphi(.,1) .\eqno(5.11)
$$ Indeed, we verify directly, that if $\tilde r = r + x e_{-1} \wedge e_1$, then (5.3) holds with 
$Y = x(1+x)K$. 

This example is exceptional. Drin'feld's Map $\phi$ is just the modifying term in (5.3), and it is not zero. But its
differential is zero, so the   Lie ``quasi bialgebra" is actually coassociative, and hence a Lie bialgebra.

(3) Let us turn to the elliptic deformation of $\widehat{sl_q(2)}$. In this case $\hat \Gamma_1 = \hat\Gamma_2$
consists of the simple roots $f_1$ and $f_0$,  
$$
\tau f_1 = f_0,\quad \tau f_0 =  f_1. $$
\vskip-5mm
\no and
$$ F = 1-\epsilon(f_1\otimes f_{-0} + f_0 \otimes f_{-1}) + o(\epsilon^2).\eqno(5.12)
$$ To use (5.4) we temporarily pass to the loop algebra, setting $c=0$ to obtain
$$
 X^m = {(-)^m\over 1-x\epsilon^{-2m}}\biggl( \eta^{-1} H \otimes H +\bigg\{ \matrix {x\epsilon^{-m}(f_1 \otimes
f_{-1} + f_0 \otimes f_{-0}),\quad m~ {\rm even},\cr
 \cr -x \epsilon^{-m}(f_1 \otimes f_{-0} + f_0 \otimes f_{-1}),\quad m ~{\rm odd}}\biggr),\eqno(5.13)
$$ with 
$$ H = (\varphi + \varphi)(1),\quad x = \lambda/\mu,\quad \eta = 2\varphi(1,1) = H(1).
$$  For the generators of the affine Lie algebra we have written
$$ 1\otimes H^n = \lambda^n (1 \otimes H), \quad H^n \otimes 1 = \mu^n (1 \otimes H),\eqno(5.14)
$$ and so on. In terms of Lie generators,
$$\eqalign{ X^m = -\hskip-1mm\sum_{n=1}^\infty\biggl( \eta^{-1}&(-\epsilon^{2n})^m ~H^n \oplus H^{-n}\cr &+
\epsilon^{m(2n-1)}\bigg\{\matrix{f_1^{n-1} \otimes f_{-1}^{1-n} + f_{-1}^n \otimes f_1^{-n},\hskip1mmm~{\rm even},\cr
\cr f_1^{n-1} \otimes f_1^{-n} + f_{-1}^n \otimes   f_{-1}^{1-n},\hskip1mm m~{\rm odd}\cr}
\biggr). 
\cr}\eqno(5.15)
$$  In this form the result holds on the full Kac-Moody algebra. 

The modification of the Yang-Baxter equation is given by (5.3) and (5.6).
 For simplicity, take the parameter $u$ to be $1/2$. Then the complete result is (with $N= 2)$
$$
\delta r \wedge Y = {\epsilon \over N} {d\over d\epsilon}  \tilde r \wedge c.\eqno(5.16)
$$ As for the deformed r-matrix we have $\tilde r = r + \delta r$, with $\delta r = X - X^t,~ X = \sum X^m$. The
sums represent elliptic functions,
\def\sn{{\rm sn}}
\def\cn{{\rm cn}}
\def\dn{{\rm dn}} 
$$\eqalign{ &\quad\quad  \tilde r  =   u c\otimes d + (1-u)d \otimes c + {i\over 2\pi}[i{\cn (u)\over \sn
(u)}(\varphi+\varphi^t)  \cr &+ {1\over \sn (u)}(e^{iu} f_1 \otimes f_{-1} + e^{-iu}f_{-1} \otimes f_1)  + {\sn
(u)\over \dn (u)}(e^{iu} f_1 \otimes f_{-0} + e^{-iu} f_0 \otimes f_{-1})]. 
\cr}\eqno(5.17)
$$
\b (4) Next, we present the elliptic quantum groups, in the sense of [Ba] and [Be].
  Let ${\cal A}'$ be the full, quantized, affine Kac-Moody algebra associated with $sl(N)$. Let the simple roots
$e_1,...e_{N-1}$ be ordered as in the Dynkin diagram, and let
$e_0$ be the imaginary root. Let $\tau e_i = e_{i+1}$ for $i = 1,...,N-2,~\tau e_{N-1} = e_0$ and $
\tau e_0 = e_1$. Suppose that $Q$ reduces to 1 on the quotient defined by the center. Then for a pair $(\s,\rho) =
(f_\s,f_{-\tau
\s})$ we have  
$$ Q(\s,\rho) = e^{-c\bigl( ud(\rho)+(1-u)d(\s)\bigr)},\quad d(\s) = \delta_\s^0.
$$
    The solution of (5.4)  gives us the  classical r-matrix $\tilde r = r + X - X^t,~ X := \sum  X^m$. It agrees
with that obtained by Belavin and Drinfeld [B][BD], except that it is augmented by the contribution of the extension,
 $$
\tilde r =  r({\rm Belavin}) + u c\otimes d + (1-u)d \otimes c.
\eqno(5.18)
$$ It satisfies the modified Yang Baxter equation (5.3), with $\delta r \wedge Y$ as in (5.16,  for all $N$.

 \vskip.5cm {\bf Remark 5.1.} To interpret (5.18) correctly, one must express it in terms of Lie generators. This is
done by   expanding the elliptic functions in terms of powers of $x =
\lambda/\mu$. It has  been emphasized that both positive and negative powers of $x$ occur in the elliptic r-matrix.
Here the  point of view of deformation theory comes to the rescue; to each   order in
$\epsilon$ there is only a finite number of negative powers of $x$. The r-matrix is a power series in $\epsilon$ and
each term is a Lorent series in $x$ with finite principal part.  
\b (5) An affine version of Example (1) is the deformation of $\widehat {sl_q(2)}$ with $\tau = 1$ and
$$ F = 1-\epsilon(kf_1 \otimes f_{-1} + f_0 \otimes f_{-0}) + o(\epsilon^2).\eqno(5.19)
$$
 The recursion relation (5.4) gives us, for the loop algebra,
$$ X^m = \Psi^mC_0 + A^me_1 \otimes  e_{-1} + B^m e_{-1} \otimes e_1,
$$ with
$$
\Psi^m = {-x^{-1}(k\epsilon^2)^m\over
 1- x^{-1}(k\epsilon^2)^m},\quad A^m = {-(k\epsilon)^m\over 1- x^{-1}(k\epsilon^2)^m} = xk^mB^m, \quad x =
\lambda/\mu,
$$ and $C_0 =\varphi + \varphi^t$. To pass to the full Kac-Moody algebra one should express this in the manner of
(5.15); note that $C_0$ does \underbar {not} include the extension. Summing over $m$ one gets
$$ X = \sum_{m=1}^\infty X^m = \Psi C_0 + A\,e_1 \otimes e_{-1} + B\, e_{-1} \otimes e_1,
$$ with
$$ 
\Psi = - \sum {(k\epsilon^2)^n\over 1 - (k\epsilon^2)^n}x^{-n},~~~  A =   - x\sum {k\epsilon(k\epsilon^2)^n\over 1 -
k\epsilon(k\epsilon^2)^n}x^{-n}, ~~~ B =  - \sum {\epsilon(k\epsilon^2)^n\over 1 - \epsilon(k\epsilon^2)^n}x^{-n}. 
$$ Finally,
$$
\tilde r = {\varphi - \varphi^t\over 2} + {1\over 2i}\,\zeta(u)C_0 + {i\over 2} Z(u)\bigl( e^{i\pi u} e_1 \otimes
e_{-1} + e^{-i\pi u} e_{-1} \otimes e_1\bigr),\eqno(5.20)
$$ with
$$
\eqalign{
\zeta(u) &= {1\over \tan \pi u} + \sum{(k\epsilon^2)^n \over 1 -  (k\epsilon^2)^n} \sin \pi u,\cr Z(u) &= {1\over
\sin \pi u} -2i \sum \biggl\{{\epsilon(k\epsilon^2)^n\over 1 - \epsilon(k\epsilon^2)^n} e^{(2n-1)\pi iu} -
{k\epsilon(k\epsilon^2)^n\over 1 - k\epsilon(k\epsilon^2)^n} e^{-(2n-1)\pi iu}\biggr\}\cr }\eqno(5.21)
$$ The modified Yang-Baxter relation takes the form (5.3) with
$$
\delta r \wedge Y = \epsilon  {d\over d\epsilon}\tilde r \wedge ({c\over 2}+K), \quad K = (\varphi +
\varphi^t)(1).
$$ The Hopf algebra is benevolent but not central. The r-matrix is similar to the one constructed in [ABB] for the
Calogero-Moser model.
\b
\no {\bf Remark 5.2.} The function here denoted $\zeta$ differs from the elliptic $\zeta$-function by the omission
of the term linear in $u \propto \log(\lambda/\mu)$. Such a term is allowed in the sense that it would make no
contribution to the modified Yang-Baxter relation, but it would destroy the connection to Kac-Moody algebras. 
 
 \b (6) Finally, we generalize the last example to an arbitrary quantized, untwisted,   affine Kac-Moody algebra ~
$\hat {}$ \hskip-2.35mm \g~ over a Lie algebra \g. Let   
$\Gamma_1 = \Gamma_2$ be the set of simple roots, take $\tau=1$ and
$$ F = 1 -  \sum_{\s \in \Gamma_1}\bigl(\epsilon_\s   f_\sigma \otimes f_{-\s} + \epsilon f_0 \otimes f_{-0}) +
o(\hbar^2),\eqno(5.22)
$$ with $\epsilon_\s = k_\s\epsilon, ~k_\s$ a set of fixed parameters. We find 
$$ X^m = \Psi^mC_0 + \sum_iA_i^mE_i \otimes  E_{-i} + B_i^m E_{-i} \otimes E_i~,
$$ with
$$
\Psi^m = {-x^{-1}(\epsilon\epsilon_+)^m\over
 1- x^{-1}(\epsilon\epsilon_+)^m},
\quad
 A_i^m = {-\epsilon_i^m\over 1- x^{-1}(\epsilon\epsilon_+)m^{-1}},\quad  B_i^m ={ -x^{-1}(\epsilon
\epsilon_+/\epsilon_i)^m\over 1- x^{-1}(k\epsilon^2)m^{-1}}.
$$ Here $\epsilon_i = \prod\epsilon_\s$, where the products runs over the simple roots that make up the root $E_i$,
and $E_+$ is the greatest root. 
 
 The modified r-matrix is
$$
\tilde r = {\varphi - \varphi^t\over 2} + {1\over 2i}\zeta(u)C_0 + 
\sum_iZ_i(u)\bigl(e^{iu}E_i \otimes E_{-i} + e^{-iu}E_{-i}\otimes E_i\bigr),\eqno(5.23)
$$ with
$$
\eqalign{
\zeta(u) &= {1\over \tan \pi u} + \sum{(\epsilon_+\epsilon_0)^n \over 1 -  (\epsilon_+\epsilon_0)^n} \sin \pi u,\cr
Z_i(u) &= {1\over \sin \pi u} -2i \sum \biggl\{{\epsilon_i(\epsilon_+\epsilon_0)^n\over 1 -
\epsilon_i(\epsilon_+\epsilon_0)^n} e^{(2n-1)\pi iu} -
{(\epsilon_+\epsilon_0/\epsilon_i)(\epsilon_+\epsilon_0)^n\over 1 -
(\epsilon_+\epsilon_0/\epsilon_i)(\epsilon_+\epsilon_0)^n} e^{-(2n-1)\pi iu}\biggr\}\cr }\eqno(5.24)
$$
 
 The extra term in the modified Yang-Baxter relation takes the form
$$
\delta r \wedge Y = \sum_\s \epsilon_\s {d\over d \epsilon_\s}  \tilde r \wedge ({c\over N} + K(\sigma)),
\quad K(\s) = (\varphi + \varphi^t)(\s),\eqno(5.25)
$$  where $N-1$ is the length of the longest root. There is no difficulty in dealing with the twisted,  affine
Kac-Moody algebras.
\b
 
Suppose that $\epsilon_0 = 0$, then $\Psi_m = B^m = 0, A_i^m = - \epsilon_i^m$ and $A_i =
-\epsilon_i(1-\epsilon_i)^{-1}$. In particular, for $sl(N)$,  
$$
\tilde r = r + \sum_{i<j}{q_i\over q_i-q_j}e_{ij} \wedge e_{ji}.
$$ This is similar, but not identical to the r-matrix in [ABB]. Note that here $\tilde r + \tilde r^t$ is the
Killing form, while the r-matrix in [ABB] is unitary, $r + r^t = 0$.
\b
\no {\bf Final Remarks.} (1) It should be emphasized that the elliptic quantum groups associated with $sl(N)$, ``in
the sense of Baxter and Belavin", form a class quite apart from the ``elliptic quantum groups" in the sense of
Felder.  The homomorphism $\tau$ (see Proposition 4.2) is in the first case a nontrivial diagram automorphism of
$sl(N)$, in the second case the identity.  (2) Unfortunately, I have been unable  to identify the functions defined
by the series in (5.21), and I do not know the genus of the curve on which they live.

\ve

\no {\bf References}
 \noindent
\parindent0pt 
 
[ABB] J. Avan, O. Babelon,   and E. Billey, The Gervais-Neveu-Felder equation and 
   \line{\hfil the quantum Calogero-Moser system. hep-th/9505091.}

[BBB] O. Babelon, D. Bernard and E. Billey, A Quasi-Hopf interpretation of quantum \line {\hfil 3-j and 6-j symbols
and difference equations, q-alg 9511019.}

[Ba]\hskip5mm R.J. Baxter, Partition Function of the eight-Vertex Model,  Annals of Physics \line {\hfil  {\bf 70} 
(1972) 193-228.}

[Be]\hskip4.5mm A.A. Belavin, Dynamical Symmetry of Integrable Systems, Nucl. Phys. {\bf 180}
  \line {\hfil  (1981)  198-200.}

[BD]\hskip2mm  A.A. Belavin and V.G. Drinfeld, Triangle Equation and Simple Lie Algebras,  Sov. \line {\hfil
Sci.Rev.Math.{\bf 4} (1984) 93-165.}

 [CG] \hskip1mm  E. Cremmer and J.-L. Gervais, Commun.Math.Phys. {\bf 134} (1990).

[D1]\hskip3mm  V.G. Drinfeld, Quantum Groups, in Proceedings of the International Congress of
   \line {\hfil Mathematicians,~Berkeley,  A.M. Gleason, ed. (A.
    M.  S.,  Providence  1987).}

[D2]\hskip3mm  V.G. Drinfeld, Quasi Hopf Algebras, Leningrad Math.J. {\bf 1} (1990) 1419-1457.

[Fe] G. Felder, Elliptic Quantum Groups, hep-th/9412207.

[Fo] O. Foda, K. Iohara, M. Jimbo, T. Miwa and H. Yan, An Elliptic Algebra for $sl(2)$, RIMS preprint 974.

[FR]~~\hskip .7mm I.B. Frenkel and N.Yu. Reshetikhin, Quantum Affine Algebras and Holonomic 
   \line {\hfil Difference Equations,Commun.Math.Phys. (1992) {\bf 146} 1-60).}

[Fr] C. Fr\o nsdal, Generalization and Deformations of Quantum Groups, to appear in \line {\hfil RIMS Publications.}
 
[FG1] C. Fr\o nsdal and A. Galindo, Contemporary Mathematics, {\bf 175} (1994) 73-88.
 
[G] \hskip3.5mm M. Gerstenhaber, Ann.Math {\bf 78} (1963) 267-288.

[GN]\hskip1mm J.-L. Gervais and A. Neveu, Novel Triangle Relation and Absence of Tachyons in 
\line {\hfil Liouville String Field Theory, Nucl.Phys. {\bf 238} (1984) 125.}

[Au] \hskip3mm H. Au-Yang, B. McCoy, J.H.G. Perk and M.-I. Yan, Phys.Lett. {\bf A123} (1987) 
\line {\hfil 219-223. B.McCoy, J.H.H. perk and C.-H. Sah, Phys.lett. {A125}(1987) 9-14. }

\end
  We illustrate the calculation of $X^m
\otimes Y^m$ by obtaining the coefficient of
$H^n
\otimes H^{-n}$, in the case of even $m$. The generator $H^n$ has the same weight as $(f_0f_1)^n$. For even $m$ we 
have $\tau^m f_1 = f_1$. The contribution to $Y^m$ associated with each pair $f_\s \otimes f_{-\s}$ is
$(m/2)[\varphi(0,.) + \varphi(.,1) + \varphi(1,.) + \varphi(.,0)] = (m/2)c$. The coveted  coefficient is thus $mnc$.